\documentclass[12pt]{article}

\usepackage{amsmath,amssymb}

\begin{document}

\title{Differential operators on equivariant vector bundles over symmetric spaces}
\author{Anton Deitmar}

\date{}
\maketitle

\tableofcontents

\def \1{{\bf 1}}
\def \a{{{\mathfrak a}}}
\def \ad{{\rm ad}}
\def \al{\alpha}
\def \ar{{\alpha_r}}
\def \A{{\mathbb A}}
\def \Ad{{\rm Ad}}
\def \Aut{{\rm Aut}}
\def \b{{{\mathfrak b}}}
\def \bs{\backslash}
\def \B{{\cal B}}
\def \c{{\mathfrak c}}
\def \cent{{\rm cent}}
\def \C{{\mathbb C}}
\def \CA{{\cal A}}
\def \CB{{\cal B}}
\def \CC{{\cal C}}
\def \CD{{\cal D}}
\def \CE{{\cal E}}
\def \CF{{\cal F}}
\def \CG{{\cal G}}
\def \CH{{\cal H}}
\def \CHC{{\cal HC}}
\def \CL{{\cal L}}
\def \CM{{\cal M}}
\def \CN{{\cal N}}
\def \CP{{\cal P}}
\def \CQ{{\cal Q}}
\def \CO{{\cal O}}
\def \CS{{\cal S}}
\def \CT{{\cal T}}
\def \CV{{\cal V}}
\def \det{{\rm det}}
\def \diag{{\rm diag}}
\def \dist{{\rm dist}}
\def \Eig{{\rm Eig}}
\def \End{{\rm End}}
\def \Fx{{\mathfrak x}}
\def \FX{{\mathfrak X}}
\def \g{{{\mathfrak g}}}
\def \ga{\gamma}
\def \Ga{\Gamma}
\def \h{{{\mathfrak h}}}
\def \Hom{{\rm Hom}}
\def \im{{\rm im}}
\def \Im{{\rm Im}}
\def \Ind{{\rm Ind}}
\def \k{{{\mathfrak k}}}
\def \K{{\cal K}}
\def \l{{\mathfrak l}}
\def \la{\lambda}
\def \lap{\triangle}
\def \La{\Lambda}
\def \m{{{\mathfrak m}}}
\def \mod{{\rm mod}}
\def \n{{{\mathfrak n}}}
\def \name{\bf}
\def \N{\mathbb N}
\def \o{{\mathfrak o}}
\def \ord{{\rm ord}}
\def \O{{\cal O}}
\def \p{{{\mathfrak p}}}
\def \ph{\varphi}
\def \prf{\noindent{\bf Proof: }}
\def \Per{{\rm Per}}
\def \q{{\mathfrak q}}
\def \qed{$ $\newline $\frac{}{}$\hfill {\rm Q.E.D.}\vspace{15pt}\pagebreak[0]}
\def \Q{\mathbb Q}
\def \res{{\rm res}}
\def \R{{\mathbb R}}
\def \Re{{\rm Re \hspace{1pt}}}
\def \r{{\mathfrak r}}
\def \ra{\rightarrow}
\def \rank{{\rm rank}}
\def \s{{\mathfrak s}}
\def \supp{{\rm supp}}
\def \Spin{{\rm Spin}}
\def \t{{{\mathfrak t}}}
\def \T{{\mathbb T}}
\def \tr{{\hspace{1pt}\rm tr\hspace{2pt}}}
\def \vol{{\rm vol}}
\def \z{\zeta}
\def \Z{\mathbb Z}
\def \={\ =\ }

\newcommand{\rez}[1]{\frac{1}{#1}}
\newcommand{\der}[1]{\frac{\partial}{\partial #1}}
\newcommand{\norm}[1]{\parallel #1 \parallel}
\renewcommand{\matrix}[4]{\left( \begin{array}{cc}#1 & #2 \\ #3 & #4 \end{array}
            \right)}

\newcounter{lemma}
\newcounter{corollary}
\newcounter{proposition}
\newcounter{theorem}

\newtheorem{conjecture}{\stepcounter{lemma} \stepcounter{corollary}
    \stepcounter{proposition}\stepcounter{theorem}
    \stepcounter{subsection}\hskip-12pt Conjecture}[section]
\newtheorem{lemma}{\stepcounter{conjecture}\stepcounter{corollary}
    \stepcounter{proposition}\stepcounter{theorem}
    \stepcounter{subsection}\hskip-7pt Lemma}[section]
\newtheorem{corollary}{\stepcounter{conjecture}\stepcounter{lemma}
    \stepcounter{proposition}\stepcounter{theorem}
    \stepcounter{subsection}\hskip-7pt Corollary}[section]
\newtheorem{proposition}{\stepcounter{conjecture}\stepcounter{lemma}
    \stepcounter{corollary}\stepcounter{theorem}
    \stepcounter{subsection}\hskip-7pt Proposition}[section]
\newtheorem{theorem}{\stepcounter{conjecture} \stepcounter{lemma}
    \stepcounter{corollary}\stepcounter{proposition}
    \stepcounter{subsection}\hskip-11pt Theorem}[section]

$$ $$

Generalizing the algebra of motion-invariant differential operators on a symmetric space we study invariant operators on equivariant vector bundles. We show that the eigenequation is equivalent to the corresponding eigenequation with respect to the larger algebra of all invariant operators. We compute the possible eigencharacters and show that for invariant integral operators the eigencharacter is given by the Abel transform. We show that sufficiently regular operators are surjective, i.e. that equations of the form $Df=u$ are solvable for all $u$.

\newpage

\section{Equivariant vector bundles} \label{equiv}
Let $X$ denote a manifold with a smooth action of a Lie group $G$. 
An {\it equivariant vector bundle} $E$ over $X$ is a smooth vector bundle $\pi : E\ra X$ together with a smooth action of $G$ on $E$ such that $\pi(gv)=g\pi(v)$ for $v\in E$ and such that all maps of the fibres $g: E_x\ra E_{gx}$, $g\in G$, $x\in X$, are linear. 
An example is given by $E=X\times V$, where $(\sigma ,V)$ is a finite dimensional representation of $G$ and $G$ acts on $E$ by $g(x,v)=(gx,\sigma(g)v)$.

Now assume $X$ to be a homogeneous space, i.e. $X=G/H$ for a closed subgroup $H$ of $G$. Given an equivariant bundle $E$ over $X$ we get a representation of $H$ on the fibre over $eH$, where $e$ denotes the neutral element of $G$. Conversely given a representation $\tau$ of $H$ on a finite dimensional vector space $V$ we let $H$ act from the right on $G\times V$ by $(g,v)h= (gh,\tau(h)^{-1}v)$ and define $E_\tau= (G\times V)/H$. This is a vector bundle over $X$. This construction gives an equivalence of categories between the category of equivariant vector bundles over $X$ and the category of finite dimensional representations of $H$.
If $\tau$ splits as a direct sum $\tau=\tau_1\oplus\tau_2$ then $E_\tau\cong E_{\tau_1}\oplus E_{\tau_2}$ and vice versa.

Now let $G$ denote a connected semisimple Lie group with finite center and let $K$ denote a maximal compact subgroup. The quotient $X=G/K$ is diffeomorphic to $\R^n$ for some $n$. Using the Killing form of $G$ one defines a Riemannian metric on $X$ such that the group $G$ acts by isometries. This is the most general symmetric space without compact factors. See \cite{helg} for further details.

\section{Differential operators}
Let $\hat{K}$ denote the set of isomorphism classes of irreducible unitary representations of the group $K$. Since $K$ is compact, every $\tau\in\hat{K}$ is finite dimensional.
We do not distinguish between a class in $\hat{K}$ and a representative.
For $(\tau ,V_\tau)\in\hat{K}$ let $E_\tau\ra G/K$ denote the vector bundle as in section \ref{equiv}. The group $G$ acts on the space of smooth sections $\Ga^\infty(E_\tau)$ of the bundle $E_\tau$ by $g.s(x)=gs(g^{-1}x)$, $x\in X$, $g\in G$, $s\in \Ga^\infty(E_\tau)$.
Hence $G$ also acts on differential operators by conjugation.
Let $\CD_\tau$ denote the algebra of $G$-invariant differential operators on $E_\tau$, i.e. those operators $D$ that satisfy $D(g.s)=g.D(s)$ for any $s\in\Ga^\infty(E_\tau)$.

Let $C_\tau^\infty(G)$ denote the space of all infinitely often differentiable maps from $G$ to $V$ with $f(kx)=\tau(k)f(x)$ for $k\in K$ and $x\in G$. The group $G$ acts on the space $C_\tau^\infty(G)$ by translations from the right. The map $s\mapsto f_s$ with $f_s(x)=xf_s(x^{-1})$ gives a $G$-isomorphism of $\Ga^\infty(E_\tau)$ to $C_\tau^\infty(G)$. We conclude that the algebra $\CD_\tau$ acts on $C_\tau^\infty(G)$.

Let the compact group $K$ act on the space $C^\infty(G)$ of smooth functions on $G$ by
$$
L_k(x)\= f(k^{-1}x),\ \ \ k\in K, \ x\in G.
$$
Then there is a decomposition into $K$-isotypes:
$$
C^\infty(G)\= \bigoplus_{\tau\in\hat{K}} C^\infty(G)(\tau).
$$
The sum here means that finite sums on the right hand side are dense in the Fr\'echet space $C^\infty(G)$. Further $C^\infty(G)(\tau)$ is the space of functions in $C^\infty(G)$ that transform under $L$ according to $\tau$.

For $(\tau, V)\in \Hat{K}$ fix any nonzero $v^*$ in the dual space $V^*$ and for $f\in C_\tau^\infty(G)$ let $Bf(x)=v^*(f(x))$. Then $Bf\in C^\infty(G)(\tau)$ and the map $B$ is a $G$-isomorphism.
We have
$$
\Ga^\infty(E_\tau)\ \cong\ C_\tau^\infty(G)\ \cong \ C^\infty(G)(\tau).
$$
Let $\g$ denote the Lie algebra of $G$. The universal enveloping algebra $U(\g)$ may be viewed as the algebra of all right invariant differential operators on $G$, i.e. the algebra of all differential operators $D$ on $G$ such that $R_gD=DR_g$ for any $g\in G$. Here for a function $\ph\in C^\infty(G)$ and $g\in G$ we have $R_g\ph(x)=\ph(xg)$ for all $x\in G$.
Let $U(\g)^K$ denote the subalgebra of all differential operators which are $K$-invariant on the left side, i.e. which satisfy $L_kD=DL_k$ for any $k\in K$.
The algebra $U(\g)^K$ leaves stable the decomposition of $C^\infty(G)$ and thus acts on $C^\infty(G)(\tau)$. 
We therefore get an algebra homomorphism
$$
\ph_\tau : U(\g)^K\ra \CD_\tau.
$$

\begin{proposition}
The homomorphism $\ph_\tau$ is surjective and the intersection of all kernels $\ker(\ph_\tau)$ for varying $\tau$ is zero.
\end{proposition}

\prf
Let the algebra $U(\g)\otimes \End(V)$ acts on the space $C^\infty(G,V)$ of smooth functions on $G$ with values in $V$. 
The group $K$ acts on $U(\g)$ via the adjoint representation and on $\End(V)$ by means of conjugation via $\tau$.
Then the algebra of $K$-invariants, $(U(\g)\otimes\End(V))^K$ acts on $C_\tau^\infty$.
The annihilator $I$ of $C_\tau^\infty$ in $U(\g)\otimes\End(V)$ is generated by elements of the form $XY\otimes T+X\otimes T\tau(Y)$ with $X\in U(\g)$, $Y\in \k$ and $T\in\End(V)$.
Since $K$ is reductive we have
$$
\CD_\tau\ \cong\ (U(\g)\otimes\End(V))/I)^K\ \cong\ (U(\g)\otimes\End(V))^K/I^K
$$
Since $\tau$ is irreducible and $K$ is connected it follows $\tau(U(k))=\End(V)$.
This implies that any element of $\CD_\tau$ can be written in the form $Z\otimes 1$ for some $Z\in U(\g)$.
It follows that $Z$ must be in $U(\g)^K$, which implies the surjectivity of $\ph_\tau$.

For the second assertion assume $X$ is in the intersection of all kernels.
Then, since $C^\infty(G) = \bigoplus_\tau C^\infty(G)(\tau)$, we get that $Xf=0$ for every $f\in C^\infty(G)$, this gives $X=0$.
\qed

\begin{corollary}
$\CD_\tau$ is finitely generated as $\C$-algebra.
\end{corollary}

\prf
$U(\g)$ has a natural filtration by order. The associated graded version equals the symmetric algebra $S(\g)$ over $\g$.
Since the adjoint action of $K$ preserves the filtration, we have a filtration on $U(\g)^K$ with graded version $S(\g)^K$.
The latter is finitely generated by invariant theory, hence the former is, too.
\qed

\noindent
{\bf Examples.}

For $\tau =1$ the trivial representation, the algebra $\CD_\tau$ is the algebra of $G$-invariant differential operators on $G/K$.
In this case $\CD_\tau$ is isomorphic to the polynomial ring in $r$ generators, where $r$ is the rank of the symmetric space $G/K$.

Let $G=SO(n,1)^+$ the group of motions on the hyperbolic space $H_n$.
For $\tau=\wedge^p(\Ad)$ the space $\Ga^\infty(E_\tau)$ is just the space of $p$-differential forms on $H_n$.
For $p=0$ or $p=n$ the algebra $\CD_\tau$ is the polynomial ring in one variable, generated by $\delta d$ and $d\delta$ respectively, where $d$ is the exterior differential and $\delta$ its formal adjoint.
For $n$ odd and $p=(n\pm 1)/2$ the algebra $\CD_\tau$ is generated by $d\delta$ and $*d$ resp. $\delta d$ and $d*$, where $*$ is the Hodge operator, with the generating relations
$$
d\delta *d = 0 = *dd\delta,\ \ \ \ \delta dd* = d*\delta d=0.
$$
In all other cases $\CD_\tau$ is generated by $d\delta$ and $\delta d$  qith the relations $d\delta \delta d=0=\delta dd\delta$.
Summarizing we get the structure of $\CD_\tau$ as:
$$
\CD_\tau\ \cong\ \left\{ \begin{array}{cl} \C[x,y]/xy & {\rm for}\ 1\le p\le n-1,\\
                                           \C[x]      & {\m for}\ p=0,n, 
                         \end{array} \right.
$$

\section{Integral operators}
A smooth function $\Phi : G\ra \End(V)$ which satisfies $\Phi(kxl)=\tau(k)\Phi(x)\tau(l)$ for $x\in G$, $k,l\in K$ is called {\it $\tau$-sherical}. The algebra $\CD_\tau$ acts on the set of $\tau$-spherical functions. Compactly supported $\tau$-spherical functions form an algebra under convolution:
$$
\Phi * \Psi (x)\= \int_G \Phi(y)\Psi(y^{-1}x)dy,
$$
where $dy$ denotes a Haar measure on $G$. This algebra is denoted by $\CA_\tau$. The algebra $\CA_\tau$ acts on $C_\tau^\infty(G)$ by $L_\Phi f= \Phi * f$ for $\Phi\in\CA_\tau$ and $f\in C_\tau^\infty$. The algebra $\CA_\tau$ contains an approximate identity.

Let $D_1,\dots ,D_n$ be a set of generators of the algebra $\CD_\tau$. For $z\in \C^n$ let 
$$
\Eig(z) \= \{ f\in C^\infty(G, \End(V)) |\ f\ {\rm is}\ \tau-{\rm spherical\ and}\ D_jf=z_jf\}.
$$

\begin{proposition}
For any $z\in\C^n$ we have $\dim\Eig(z)\le 1$.
\end{proposition}

\prf
By Lemma 1 in \cite{har1} any $f\in\Eig(z)$ is analytic. 
So let $f\in\Eig(z)$ and $H\in\g$.
If $H$ is small enough then
$$
f(\exp(H))\= \sum_{n\ge 0} \rez{n!}\left( \frac{\partial}{\partial t}\right)^n f(\exp(tH))|_{t=0} \= \sum_{n\ge 0} \rez{n!}H^nf(0).
$$
We get
$$
\int_K \tau(k)f(\exp(H))\tau(k^{-1})dk\= \sum_{n\ge 0}\int_K (\Ad(k)H)^nf(e) dk.
$$
Now $\int_K(\Ad(k)H)^n dk$ is a right invariant differential operator mapping $C_\tau(G)$ to itself, hence defines an element of $\CD_\tau$.
So the values $ \int_K(\Ad(k)H)^nf(e) dk$ only depend on $z$ and not on $f$.
We conclude that the function $\tr f(x)$ is determined by $z$ up to scalar. But the map $\CA_\tau\ra C^\infty(G)$, $f\ra\tr(f)$ is injective \cite{wa2}.
\qed

\begin{theorem}
Suppose $f\in C^\infty_\tau(G)$ is an eigenform for any $D\in\CD_\tau$. Then $f$ is an eigenform for every $T\in\CA_\tau$ with an eigenvalue only depending on $T$ and the eigenvalues on $\CD_\tau$.
\end{theorem}

\prf
For $\ph\in C^\infty_\tau(G)$ and $w\in V$ define $\ph_w\in C^\infty(G,\End(V))$ by
$$
\ph_w(x)v\= \langle v,w\rangle \ph(x),
$$
where $\langle .,.\rangle$ denotes the scalar product on $V$.
Fot $psi : G\ra \End(V)$ let 
$$
M\psi(x)\= \int_K\psi(xk)\tau(k^{-1})dk.
$$
Now let $f$ be as in the theorem and assume $f(e)\ne 0$. (Otherwiese replace $f$ by $R_gf$ for a suitable $g\in G$.)
Fix some $w\in V$ such that $\tr f_w(e)=1$.
Now $M(f_w)$ lies in $\Eig(z)$ for some $z$.
Let $\Phi$ be $\tau$-spherical and compactly supported.
Since $L_\Phi M(f_w)=M((L_\Phi f)_w)$ we see that there is a $\la\in\C$ such that
$$
M((\la f-L_\Phi f)_w)\= 0,
$$
and $\la$ does not depend on $f$ or $w$.
The claim now follows by the proposition.
\qed

Let $\a$ denote the Lie algebra of the maximal $\R$-split torus $A$ of $G$.
Let $G=KNA$ be a corresponding Iwasawa decomposition. 
For $\la\in\a^*_\C$ the complex dual space of $\a$ and any $v\in V$ let
$$
p_{\la,v}(kna) \= \tau(k)e^{\la(\log(a))} v.
$$
This defines $p_{\la ,v}\in C_\tau^\infty(G)$.
Let $P=MAN$ be the minimal parabolic given by $A$ and $N$. Then $M$ is a closed subgroup of $K$.

\begin{lemma}\label{3.5}
For every $T\in \CA_\tau\oplus\CD_\tau$ and every $\la\in\a_\C^*$ there is a $S_\la(T)$ in $\End(V)$ such that for all $v\in V$
$$
T(p_{\la ,v})=p_{\la ,S_\la(T)v}.
$$
We have $S_\la(T)\tau(m)=\tau(m)S_\la(T)$ for all $m\in M$
\end{lemma}

\prf
The lemma is clear by group invariance and the fact that $A$ normalizes $N$.
\qed

For a simultaneous eigenform $f\in C_\tau^\infty(G)$ of $\CD_\tau$ let $\chi_f$ denote the eigencharacter $\chi_f : \CA_\tau\oplus\CD_\tau\ra \C$ defined by
$$
Tf\= \chi_f(T)f.
$$

\begin{theorem}
Let $f$ denote a bounded simultaneous eigenform.
Then there is a $v\in V$ and a $\la\in\a_\C^*$ such that $p_{\la,v}$ is an eigenform and
$$
\chi_f\= \chi_{p_{\la,v}}.
$$
\end{theorem}

\prf
The character $\chi_f$ is determined by its values on $\CA_\tau$.
Let $p$ denote the trivial seminorm on $G$, i.e. $p(g)=1$ for all $g\in G$.
Set 
$$
\norm{\Phi}_p=\int_G\norm{\Phi(y)}p(y)dy,
$$ 
where $\norm{.}$ is the norm on $\End(V)$.
Assume $\norm{f(e)}=1$ and $\norm{f(x)}\le M$, $x\in G$.
Then we get for $\Phi\in\A_\tau$:
\begin{eqnarray*}
|\chi_f(L_\Phi)| &=& \norm{\int_G \Phi(y)f(y^{-1})dy}\\
	&\le & M\int_G\norm{\Phi(y)}dy\\
	&\le& M\norm{\Phi}_p.
\end{eqnarray*}
Thus $\chi_f$ is a $p$-continuous representation  of $\CA_\tau$.
The claim now follows from the theorem of Glover \cite{wa2}, p.40.
\qed

We now give the computation of the eigencharacters of $\CA_\tau$.
Let $\rho=\rez{2}\sum_{\alpha >0}m_\alpha \alpha\in\a^*$ be the usual modular shift, i.e. the sum runs over all positive roots and $m_\alpha$ is the dimension of the root space to the root $\alpha$.

\begin{theorem}
Let $S_\la(L_\Phi)$ denote the endomorphism of Lemma \ref{3.5}. Then with
$$
g_{\Phi}(a) \= a^{-\rho}\int_N \Phi(na) dn,
$$
(Abel transform), we get
$$
S_\la(L_\Phi)\= \int_A  a^{\rho-\la} g_{\Phi}(a) da,
$$
(Fourier transform on $A$).
Moreover, $g_\Phi(a)$ is in the center of $\tau|_M$ and $g_{\Phi *\Psi}=g_{\Phi}*g_\Psi$ with $A$-convolution on the right hand side.
The map $\Phi\ra g_\Phi$ is injective.
\end{theorem}

\prf
A calculation using the integral formula of the Iwasawa decomposition gives the first claim.
The injectivity is proved in \cite{wa2}, p.35.
\qed

\section{Surjectivity of differential operators}
Let $\theta$ denote the Cartan involution fixing $K$ pointwise.
Let $\Ga_1,\dots,\Ga_r$ denote a complete system of nonconjugate $\theta$-stable Cartan subgroups of $G$ and let $A_i=\Ga_i\cap\exp(\p)$, where $\g =\k\oplus\p$ denotes the polar decomposition of $\g$.
Let $A$ denote one of the $A_i$.
Let $L=MA$ denote the centralizer of $A$ in $G$.
Let $\tau_M$ denote the restriction f $\tau$ to $K_M=K\cap M$.
Let $C^\infty(MA,\tau)$ denote the set of $\tau_M$-spherical functions on $MA$.
For $g\in C^\infty(MA,\tau_M)$ let $g^\#(kman)=\tau(k)g(ma)$.
The set of these functions $g^\#$ is stable under $\CD_\tau$ and we get a homomorphism
$$
\ga : \CD_\tau\ra\CD_{\tau_{M}}(MA)\= \CD_{\tau_M}\otimes U(\a).
$$
Now every $T\in\CD_{\tau_M}$ operates on the finite dimensional space of rapidly decreasing cusp forms $^0\CC(M,\tau_M)$ as defined in \cite{har2}.
So $\ga(D)$ defines an element $^0\ga(D)\in |End(^0\CC(M,\tau_M))\otimes U(\a)$.
Let $\det(^0\ga(D))\in U(\a)$ denote the determinant.
Call $D$ {\it regular} if $\det(^0\ga(D))\ne 0$ for all $A=A_i$.

\begin{theorem}\label{4.2}
If $D\in\CD_\tau$ is regular, then $D$ is surjective as a map from $C_\tau^\infty$ to $C_\tau^\infty$.
\end{theorem}

\prf
We are going to formulate a vector bundle version of Holmgren's uniqueness theorem.
Let $X$ denote a real analytuc manifold and $E$ an analytic vector bundle over $X$.
Let $P$ be a differential operator on $E$ with analytic coefficients.
For a susbset $A$ of $X$ let $N(A)$ denote the set of normal vectors to $A$ in $T^*X$, (see \cite{hoer}, chap. 8).
Consider the determinant of the principal symbol $\sigma_P$ as a map:
$$
\det \sigma_p : T^*X\ra \C.
$$

\begin{proposition}
(Holmgren's principle)
Let $u\in \Ga_c^\infty(E)'$ be a distribution with $Pu=0$.
Then we have
$$
\det\sigma_P(N(\supp(u)))=0.
$$
\end{proposition}

\prf
Since the assertion is local in nature the proof for the trivial bundle (\cite{hoer}, Theorem 8.6.5) carries over to the present case.
\qed

Let $P$ be as above. A point $x\in X$ is called a {\it singular point} of $P$ if
$$
\det\sigma_P(T^*_xX)\= 0.
$$
Now consider the case when $X$ is a symmetric space $G/K$ and let $E=E_\tau$ be a homogeneous bundle.
Assume that $P$ is invariant, i.e. $P\in\CD_\tau$.
Let $\nabla$ denote the canonical homogeneous connection on $E$.
For $\la\in\a_\C^*$ let $p_\la(ank)=e^{\la(\log(a))}$.

\begin{lemma}
There is a section $\ga_P$ of the bundle $S(\a)\otimes \End(E)$ such that
$$
P(p_\la s)\= \ga_P(\la)p_\la s,
$$
for all parallel sections $s$ of $E$.
For this section $\ga_P$ and any parallel $s$ we have
$$
\sigma_P(dp_\la)s\= \ga_{P,m} s,
$$
where $m=\deg(P)$ and $\ga_{P,m}$ is the principal part of $\ga_P$ with respect to the gradation of $S(\a)$.
\end{lemma}

\prf
The first part is well known, the second follows from a calculation in Iwasawa coordinates.
\qed

A differential operator $P$ on $E$ is called {\it $D$-convex}, if for any Weyl group stable compact convex subset $\s$ of $\a$ and anny section $s$ of $E$ with compact support and $\supp(Ps)\subset K\exp(\s)K$ we already have $\supp(s)\subset K\exp(\s)K$.

\begin{theorem}\label{last}
Any $P\in \CD_\tau$ with $\det\sigma_P\ne 0$ is $D$-convex.
\end{theorem}

\prf
Let $P$ be as in the theorem.
Since $\det\sigma_P$ is $G$-invariant, it follows from $\det\sigma_P\ne 0$ that $P$ has no singular points.
Let for $x\in X$, $x=k\exp(H)K$, $H\in\a$,
\begin{eqnarray*}
\delta(x) &=& \inf \{ t>0  | \frac{H}{t}\in\s \},\\
\delta(A) &=& \sup \{ \delta(x) | x\in A\},\ \ \ {\rm for}\ A\subset X.
\end{eqnarray*}
Assume $\delta(\supp Ps)=1$ and $\delta(\supp s)>>1$.
Let $x_0\in\supp s$ with $\delta(x_0)=\alpha$ and $x_0=\exp(H)$, $H$ in the positive Weyl chamber of $\a$.
Let $\la\in\a_+^*$.
With Kostants convexity theorem we get
$$
p_\la(x)\ \le\ p_\la(x_0),
$$
for all $x\in\supp s$.
So $dp_\la(x_0)\in N(\supp s)$.
Hence $\det\sigma_P(dp_\la(x_0))=0$ for all $\la>>0$, hence for all $\la$.
By group invariance this gives that $x_0$ is a singular point for $P$, a contradiction.
This proves \ref{last}.
\qed

Now regular operators $P$ satisfy the condition of \ref{last}.
They further admit fundamental solutions \cite{D}.
Now \ref{4.2} follows as in \cite{helg4}.
\qed

\noindent
{\small School of Mathematical Sciences\\
    University of Exeter\\
    Laver Building, North Park Road\\
    Exeter, EX4 4QE\\
    Devon, UK}
\end{document}